\documentstyle[11pt]{article}
\newenvironment{beweis}{{\it Proof.}\ }{\ $\ \ \ \Diamond$ \\\ }

 \newcounter{nsatz}[section]
 \newcounter{nlemma}[section]
 \newcounter{ndef}[section]
 \newcounter{nhyp}[section]
 \newcounter{nconjecture}[section]
 \newcounter{ncor}[section]
 \newcounter{nrem}[section]
 \newcounter{nexample}[section]
 \newcounter{nprop}[section]
 \newcounter{nproblem}[section]

 \newenvironment{nsatz}{\refstepcounter{nsatz}{\bf \arabic{section}.\arabic{nsatz}}\
                {\sc\bf Theorem.\ }\it}{\\\\ \rm}

 \newenvironment{nconjecture}{\setcounter{nconjecture}{\value{nsatz}}\refstepcounter{nconjecture}
                \setcounter{nsatz}{\value{nconjecture}}
                {\bf \arabic{section}.\arabic{nsatz}}\
                {\sc\bf Conjecture.\ }\it}{\\\\ \rm}

 \newenvironment{nproblem}{\setcounter{nproblem}{\value{nsatz}}
                \refstepcounter{nproblem}
                \setcounter{nsatz}{\value{nproblem}}
                {\bf \arabic{section}.\arabic{nsatz}}\
                {\sc\bf Problem.\ }}{\\\\ \rm}

\begin{document}
\newcommand{\n}{{\mbox{\rm I$\!$N}}}
\newcommand{\z}{{\mbox{{\sf Z\hspace{-0.4em}Z}}}}
\newcommand{\R}{{\mbox{\rm I$\!$R}}}
\newcommand{\Q}{{\mbox{\rm I$\!\!\!$Q}}}
\newcommand{\C}{{\mbox{\rm I$\!\!\!$C}}}

\thispagestyle{empty}
\setlength{\parindent}{0pt}
\setlength{\parskip}{5pt plus 2pt minus 1pt}

\thispagestyle{empty}
\newcommand{\Syl}{{\mbox{\rm Syl}}}
\newcommand{\cl}{{\mbox{\rm cl}}}
\newcommand{\Irr}{{\mbox{\rm Irr}}}
\newcommand{\sIrr}{{\mbox{\scriptsize\rm Irr}}}
\newcommand{\Char}{{\mbox{\rm Char}}}
\mbox{\vspace{4cm}}
\vspace{4cm}
\begin{center}
{\bf \Large\bf  A combinatorial problem arising in group theory\\}
\vspace{3cm}
by\\
\vspace{11pt}
Thomas Michael Keller\\
Department of Mathematics\\
Texas State University\\
601 University Drive\\
San Marcos, TX 78666\\
USA\\
e--mail: keller@txstate.edu\\
\vspace{1cm}
\end{center}
\thispagestyle{empty}
\newpage

\begin{center}
\parbox{12.5cm}{{\small
{\sc Abstract.}
We consider a combinatorial problem occurring naturally in a group theoretical setting and
provide a constructive solution in the smallest open case.\\
}}
\end{center}
\normalsize

\section{Introduction}\label{section0}

In this paper we want to consider a combinatorial problem whose origin is in group theory, but which can
be considered and studied in a purely combinatorial context. The task is to fill a grid with $k$ rows and
infinitely many columns with integers according to certain rules. In other words, the goal is to 
determine whether a matrix $(a_{ij})$ ($i=1,\ldots, k$; $j\in\n$) with integer entries exists which
obeys some given rules.\\ 

While the statement of the problem is quite technical,
trying to solve it turns out to be an intriguing task, particularly because intuitively
it seems almost obvious that a solution is always possible, albeit so far a general 
proof of the main conjecture (which is stated below in detail) remains elusive.\\

Before providing more background information, we now state the problem formally.\\

\begin{nproblem}\label{prob1}
Let $n\in\n$ and $k\in\n$. Suppose that we are given sets of natural numbers
$S_{ij}=S_{ij}(k)\subseteq\n$ for $i=1,\ldots,k$ and $j\in\n$ such that
$|S_{ij}|=k$ for all $i=1,\ldots,k$ and $j\in\n$. Is it possible to find
$a_{ij}=a_{ij}(k,n)\in\n$ for $i=1,\ldots,k$ and $j\in\n$ such that the following
hold.\\

(1) $a_{ij}\in S_{ij}$ for $i=1,\ldots,k$, $j\in\n$\\

(2) For any $i\in\{1,\ldots,k\}$ we have the following\\

$\quad$   (a) $a_{i1},a_{i2},\ldots,a_{in}$ are mutually distinct.\\

$\quad$   (b) $a_{ij}\not\in\{a_{i1},a_{i2},\ldots,a_{i,n-1}\}$ for all
   $j\geq n$\\

(3) For any $r_1,r_2\in\{1,\ldots,k\}$ with $r_1\not=r_2$ and any $j\in\n$ we have
\[\{a_{r_11},a_{r_12},\ldots,a_{r_1j}\}\not=\{a_{r_21},a_{r_22},\ldots,a_{r_2j}\}\]
\end{nproblem}

Being able to answer the problem in the affirmative for $n=5$ and $k=9$ plays a crucial role in a group theoretical context in
\cite{keller1} and is first explicitly discussed in \cite[Section 3]{keller2}. The solution in this case actually provides the skeleton of
an inductive process used to exhibit the existence of many orbits of different sizes in certain linear group actions, and this process
is used in several places in \cite{keller1}. A solution to the
problem for smaller values of $k$ would therefore improve several results in \cite{keller1}.
We refer the reader to \cite{keller2} for more details on the algebraic significance of this problem.\\

 From the discussion in \cite{keller2}
it follows that \ref{prob1} can always be solved in case that $k\geq 2n-1$,
whereas, on the other hand, if $k\leq n$, then one can find sets $S_{ij}$, such that \ref{prob1}
cannot be solved. For the convenience of the reader we present the easy proofs of these two claims here,
but first we want to introduce one piece of terminology.\\

Given $k,n\in\n$, we say that \ref{prob1} {\it has a general solution for $k$ and $n$}
if \ref{prob1} can be solved for any choice of the $S_{ij}$.
Otherwise, we say that \ref{prob1} does not have a general solution for $k$ and
$n$.\\

We now can prove the two claims above. First, we show that \ref{prob1} does not have a general solution for $k=n$.
(From this it is an immediate consequence that \ref{prob1} does not have a general solution whenever $k\leq n$.)
To see this, simply suppose that $S_{ij}=\{1,\dots , k\}$ for all $i=1,\ldots,k$ and $j\in\n$. Then (2) forces that
$\{a_{11},\ldots,a_{1k}\}=\{1,\dots , k\}$, and likewise (2) forces $\{a_{21},\ldots,a_{2k}\}=\{1,\dots , k\}$,
but this contradicts (3) with $r_1=1$, $r_2=2$, and $j=k$.\\

Second, we show that \ref{prob1} has a general solution for $k$ and $n$ whenever $k = 2n-1$. (From this it 
follows quickly that \ref{prob1} has a general solution whenever $k\geq 2n-1$.) To see this, observe that
without loss of generality we may assume that $a_i1=i$ for $i=1,\ldots, k$. Then we choose the $a_{ij}\in S_{ij}$
for $i=1,\ldots, k$ and $j\geq 2$ subject to the following conditions: (2) must be satisfied, and 
$a_{ij}\not\in \{i+n, i+n+1,\ldots, i+2n-2\}$, where the elements in the latter set are to be read modulo $2n-1$.
As $|S_{ij}|=2n-1$, clearly the $a_{ij}$ can indeed be chosen to satisfy these conditions, and it is then
not hard to verify that (3) holds; for example, (3) holds for Rows 1 and 2, because 1 is the first entry in Row 1,
but 1 does not occur in Row 2; and (3) holds for Rows 1 and $k=2n-1$, as $2n-1$ is the first entry in Row $k$,
but does not occur in Row 1, etc. \\ 

Therefore the open question now is: What happens when $n+1\leq k\leq 2n-2$ (and $n\geq 3$)?
Does \ref{prob1} have a general solution then?
We believe that the answer is yes, which we state as the following
conjecture.\\

\begin{nconjecture}\label{con1}
The problem has a solution whenever $k=n+1$. (Thus there is a solution whenever $k\geq n+1$.)
\end{nconjecture}
The reason for this conjecture is that in the seemingly most difficult and tightest case, namely
when all the $S_{ij}$ are equal, it is easy to show that a solution exists, as we can see as follows.
Let $k=n+1$ and suppose that the the $S_{ij}$ are all equal. Without loss of generality we may assume that
the $S_{ij}=\{1,\ldots, k\}$ for all $i, j$. Then let $a_{ij}=i+j-1$ (to be read modulo $k$) for 
$i=1,\ldots, k$ and $j=1,\ldots,k-1=n$, and let $a_{ij}=a_{in}$ for 
$i=1,\ldots, k$ and $j\geq n$. Then it is easy to check that this is indeed a solution. \\ 
So in the seemingly hardest case there is an easy solution, but a general proof of the conjecture is yet
to be found.\\

The purpose of this paper is to prove this conjecture when $n=3$ which can be extended to a general result
(see the remark following \ref{satz1}). \\

\section{The proof}\label{section2}

In this section we present a solution of the problem in the smallest, as of yet still open case, of \ref{prob1},
namely $n=3$. In this case, from the general results discussed in Section 1 we know that there does not always (i.e.,
for any choice of the sets $S_{ij}(k)$) exist a solution when $k\leq 3$, and there is always a solution
when $k\geq 5$. So the open question here is: Does there always exist a solution for $k=4$? We will show that the answer
is yes, thereby confirming \ref{con1} in this case.\\

\begin{nsatz}\label{satz1}
Let $n=3$, $k=4$ and let $S_{ij}$ ($j\in\n,\ i=1,2,3,4$) be subsets of $\n$ such that $|S_{ij}|=4$ for
all $i,j$.\\

Then there exist $a_{ij}\in\n$ ($i=1,2,3,4$; $j\in\n$) such that (1), (2), (3) of \ref{prob1}
hold.
\end{nsatz}

\begin{beweis}
First for convenience we introduce some notation.\\

The hardest rule to verify is (3), because (3) describes very many conditions. Hence we introduce some
more notation making it easier to discuss (3). For $i\in\{1,2,3,4\}$ and $j\in\n$ we define
\[M(i,j)=\{a_{i,1},a_{i,2},\ldots,a_{i,j}\}.\]
Next for any $i_1,i_2\in\{1,2,3,4\}$ and any $j\in\n$ we say that
\["P(i_1,i_2,j)\mbox{  is true}"\]
if and only if
\[M(i_1,j)\not=M(i_2,j).\]
Thus (3) is satisfied if and only if $P(i_1,i_2,j)$ is true for all $i_1,i_2\in\{1,2,3,4\}$
with $i_1\not=i_2$ and all $j\in\n$. We also say , given $i_1,i_2\in\{1,2,3,4\}$, that
\[\mbox{"}Q(i_1,i_2)\mbox{ is true"}\]
if and only if $P(i_1,i_2,j)$ is true for all $j\in\n$.\\

Hence (3) is satisfied if and only if $Q(1,2)$, $Q(1,3)$, $Q(1,4)$, $Q(2,3)$, $Q(2,4)$,
and $Q(3,4)$ are all true.\\

We show how to choose $a_{ij}\in S_{ij}$ ($i=1,2,3,4$; $j\in\n$)
in several steps.\\

Step 1: Choose $a_{11}\in S_{11}$ arbitrarily. Without loss of generality, we may
assume that $a_{11}=1$.\\

Step 2: We consider two cases:\\

Case 2a: $(S_{1j}\cap S_{i1})-\{1\}=\emptyset$ for all $i\in\{2,3,4\}$ and all
$j\geq 2$. In this case we choose $a_{1j}\in S_{1j}$ for $j\geq 2$
in such a way that (2) is satisfied, which is easily possible since we have at least
three choices for $a_{12}$ and at least two choices for each $a_{1j}$ for $j\geq 3$.
Observe that in this case $Q(1,2)$, $Q(1,3)$, and $Q(1,4)$ will be true, as
$a_{i1}\not\in M(1,j)$ for all $j\in\n$.\\

So when we continue choosing the $a_{ij}$ for $i\geq 2$ and $j\in\n$, we do not
have to pay attention to the first row of the matrix $(a_{ij})$.\\

Case 2b: If we are not in Case 2a, then we can choose $r\in\n$ minimal such that
\[(S_{1r}\cap S_{i_01})-\{1\}\not=\emptyset\]
for some $i_0\in\{2,3,4\}$. It is then no loss of generality to assume that
$i_0=2$ and that $2\in S_{1r}\cap S_{21}$.\\

Then pick $a_{1r}=2$. The remaining $a_{1j}$ for $2\leq j\leq r-1$ and
$j\geq r+1$ will be chosen later.\\

Step 3: Now we turn to the second row of the matrix $(a_{ij})$. If we are in Case 2b,
we let $a_{21}=2$. If we are in Case 2a,
we may choose $a_{21}\in S_{21}$ arbitrarily, but then without loss of generality
we may assume that $a_{21}=2$.\\

Hence in any case $a_{21}=2$.\\

Now put
\[S_{2j}^*=S_{2j}-\{1\}\]
for all $j\geq 2$..\\

Then pick $a_{2j}\in S_{2j}^*$ for $j\geq 2$ in such a way that (2) is satisfied.
This is easily possible since $|S_{2j}^*|\geq 3$, leaving at least two possibilities for
$a_{22}$ and at least one possibility for $a_{2j}$ when $j\geq 3$, as then any $a_{2j}\in S_{2j}^*-\{2,a_{22}\}$
will work to satisfy (2).\\

At this point we have filled the second row making sure that 1 does not occur in it, thus obtaining that
$Q(1,2)$ is true.\\

Step 4: We now once more have to consider two cases:\\

Case 4a: $a_{2j}\not\in S_{i1}$ for all $j\geq 2$ and
$i\in\{3,4\}$.\\

Then no matter how we later choose $a_{31}$ and $a_{41}$ (note that we must choose them different
from 1 or 2 so that $P(1,3,1)$ and $P(1,4,1)$ hold), we will always have $a_{i1}\in M(i,j)$
and $a_{i1}\not\in M(2,j)$ for $i\in\{3,4\}$ and $j\in\n$, and hence we know that $Q(2,3)$
and $Q(2,4)$ automatically hold in this case. Now we choose $a_{31}\in S_{31}-\{1,2\}$
arbitrarily. So without loss of generality we may assume that
$a_{31}=3$.\\

Case 4b: If Case 4a does not hold, then we may choose $s\in\n$ with $s\geq 2$ minimal such that there is
an $i\in\{3,4\}$ with $a_{2s}\in S_{i1}$. Recall that by the way we filled the second row of the matrix
$(a_{ij})$ we know that $a_{2s}\not\in\{1,2\}$. So now without loss of generality we may assume that
$a_{2s}=3$ and that $i=3$. Then we let
$a_{31}=3$.\\

Hence in any case we have $a_{31}=3$.\\

Step 5: We next explain how to complete the third row. Let $S_{3j}^*=S_{3j}-\{2\}$ for $j\in\n$.
Then in both Cases 4a and 4b we just pick $a_{3j}$ ($j\geq 2$) arbitrarily out of $S_{3j}^*$ in such a way that also (2) is fulfilled;
since $|S_{3j}^*|\geq 3$ for all $j$, this is clearly possible. \\

So at this point we have completed the third row in such a way that 2 does not occur in this row. Hence
while $2\in M(2,j)$ for all $j\in\n$, we have $2\not\in M(3,j)$ for all $j\in\n$, which shows that
$Q(2,3)$ holds.\\

We also claim that $Q(1,3)$ holds $(*)$.\\

To see this, recall that in Case 2a we already know that $Q(1,3)$ holds. So we may assume that
we are in Case 2b. Then note that by the choice of $r$ we know that $a_{1j}\not=3$ for $j=1,\ldots,r$
(no matter how we will choose those $a_{1j}\in S_{1j}$ later), so we have $3\in M(3,j)$ and
$3\not\in M(1,j)$ for $j=1,\ldots,r$. Thus $P(1,3,j)$ holds for $j=1,\ldots,r$.
However, as $a_{1r}=2$ (see Step 2, Case 2b), we will have $2\in M(1,j)$ for all $j\geq r$
(no matter how we choose these $a_{ij}\in S_{ij}$ later), and we made sure that $2\not\in M(3,j)$
for all $j\in\n$. Hence $P(1,3,j)$ also holds for all $j\geq r$, and thus altogether $Q(1,3)$ holds,
completing the proof of $(*)$.\\

Step 6: Let us pause for a moment and see what we have accomplished so far. We have
completed Rows 2 and 3, and we have chosen at least parts of Row 1 so that in any case we have
that $Q(1,2)$, $Q(1,3)$, and $Q(2,3)$ are
true.\\

If we are in Case 2a, then we even have completed Row 1 as well and know that $Q(1,4)$ is true, so then
it remains to fill Row 4 such that (2) is satisfied and in addition, $Q(2,4)$ and $Q(3,4)$ will hold.
(If in addition we are in Case 4a, then $Q(2,4)$ holds and we only have to make sure that
$Q(3,4)$ holds.)

If we are in Case 2b, then we have to fill Row 4 and to complete Row 1 in such a way that (2) holds
for both rows and that also $Q(1,4)$, $Q(2,4)$, and $Q(3,4)$ are true. (If in addition we are in Case 4a,
then $Q(2,4)$ holds and we only need to make sure that $Q(1,4)$ and $Q(3,4)$
hold.)\\

Step 7: As we start working on the fourth row of $(a_{ij})$, we once more consider two cases. Recall that
$a_{3j}\not\in\{2,3\}$ for all $j\geq 2$.\\

Let $a_{41}\in S_{41}-\{1,2,3\}$. Since $|S_{41}|=4$, that is possible. Note that then $P(i,4,1)$ will
be satisfied for $i=1,2,3$. Without loss of generality, we may assume that $a_{41}=4$.
Now let $S_{4j}^*=S_{4j}-\{3\}$ for $j\in\n$, and then pick $a_{4j}\in S_{4j}^*$ arbitrarily for $j\geq 2$ in
such a way that (2) is satisfied. As $|S_{4j}^*|\geq 3$, this is obviously
possible.\\

So we have completed Row 4 in such a way that (2) holds and 3 does not occur in it. Hence $3\in M(3,j)$
for all $j\in\n$ and $3\not\in M(4,j)$ for all $j\in\n$ which shows that $P(3,4,j)$ is true for $j\in\n$.
Thus $Q(3,4)$ is satisfied.\\

Next we claim that $Q(2,4)$ holds. $(**)$\\

To prove this, first assume we are in Case 2a.\\

If in addition we are in Case 4a, then we already know that $Q(2,4)$ holds. Hence we now may
assume that we are in Case 4b. Then from the definition of $s$ in Step 4, Case 4b, we know that
$a_{2j}\not=4$ for $j=1,\ldots,s-1$ and since $a_{2s}=3$, we even have $a_{2j}\not=4$ for $j=1,\ldots,s$.
Hence $4\in M(4,j)$ and $4\not\in M(2,j)$ for $j=1,\ldots,s$, so $P(2,4,j)$ holds for
$j=1,\ldots,s$. Moreover, as $a_{2s}=3$, we have $3\in M(2,j)$ for $j\geq s$, whereas we made sure
that $3\not\in M(4,j)$ for $j\in\n$. Thus $P(2,4,j)$ holds for $j\geq s$. Altogether we conclude that
$Q(2,4)$ holds and so $(**)$ is proved.\\

Step 8: We now finish the proof by finally completing Row 1, as
needed.\\

First suppose that we are in Case 2a. By what we saw in Step 6 and Step 7 we have determined the entire
matrix such that (1), (2), (3) hold so that we are done in this
case.\\

So it remains to consider Case 2b, and by Step 6 and our work done in Step 7 it remains to choose
$a_{1j}\in S_{1j}$ for $2\leq j\leq r-1$ and $j\geq r+1$ in such a way that Row 1 satisfies (2) and
$Q(1,4)$ holds.\\

To do this, let $S_{1j}^*=S_{1j}-\{4\}$. Then pick $a_{1j}\in S_{1j}^*$ for $2\leq j\leq r-1$ and
$j\geq r+1$ in such a way that (2) holds for Row 1. As $|S_{ij}^*|\geq 3$, this is clearly possible.
So we have completed Row 1 in such a way that 4 does not occur in it. Hence $4\not\in M(1,j)$
for $j\in\n$. On the other hand, $4\in M(4,j)$
for all $j\in\n$. Hence $Q(1,4)$ holds, and the proof is complete.
\end{beweis}

Remark: It is not too hard to use the ideas in the above proof to establish the new
general upper bound that \ref{prob1} has a solution whenever $k=2n-2$. This is an exercise
we leave to the interested reader.\\


\begin{thebibliography}{99}
\bibitem{keller1} T. M. Keller, Orbit sizes and character degrees, II,
J. reine angew. Math. {\bf 516} (1999), 27--116.
\bibitem{keller2} T. M. Keller, Orbits in finite group actions,
Groups St. Andrews 2001 in Oxford. Vol. II, London Math. Soc. Lecture Note Ser., 305, 
Cambridge Univ. Press, Cambridge, 2003, 306–331.


\end{thebibliography}
\end{document}